\documentclass[a4paper,11pt]{article}
\usepackage{latexsym}
\usepackage{amssymb}
\usepackage{enumerate}
\usepackage{amsfonts}
\usepackage{amsmath}
\usepackage[dvipdfmx]{graphicx} 
\usepackage[paper=a4paper,left=30mm,right=20mm,top=25mm,bottom=30mm]{geometry}
    
\newcommand{\qed}{$\Box$}

\newenvironment{@abssec}[1]{%
    \if@twocolumn

      \section*{#1}%
    \else

      \vspace{.05in}\footnotesize
      \parindent .2in
 {\upshape\bfseries #1. }\ignorespaces
    \fi}

    {\if@twocolumn\else\par\vspace{.1in}\fi}

\newenvironment{keywords}{\begin{@abssec}{\keywordsname}}{\end{@abssec}}

\newenvironment{AMS}{\begin{@abssec}{\AMSname}}{\end{@abssec}}

\newcommand\keywordsname{Key words}
\newcommand\AMSname{AMS subject classifications}
\newcommand\AMname{AMS subject classification}
\newtheorem{theorem}{Theorem}
 \newtheorem{lemma}[theorem]{Lemma}
 \newtheorem{corollary}[theorem]{Corollary}

\def\qed{\vbox{\hrule height0.6pt\hbox{%
  \vrule height1.3ex width0.6pt\hskip0.8ex
  \vrule width0.6pt}\hrule height0.6pt
 }}

\title{Interaction between initial behavior of temperature and the mean curvature of the interface in two-phase heat conductors\thanks{This research was partially supported by the Grant-in-Aid
for Scientific Research  (C) ($\sharp$ 22K03381)  of
Japan Society for the Promotion of Science. }}

\author{Shigeru Sakaguchi\thanks{Institute for Excellence in Higher Education, Tohoku University, 
Sendai, 980-8576, Japan (sigersak@tohoku.ac.jp).}}

\date{}
\begin{document}
\maketitle

\begin{abstract}
We consider the Cauchy problem for the heat diffusion equation in the whole Euclidean space consisting of two media locally with different constant conductivities, where initially one medium  has temperature 0 and the other  has temperature 1. 
Under the assumption that  a part of the interface between two media with different constant conductivities is of class $C^2$ in a neighborhood of a point  $x$ on it, we extract the mean curvature of the interface at $x$ from the initial behavior of temperature at $x$.  This result is  purely local in space. As a corollary, when the whole Euclidean space consists of two media globally with different constant conductivities, it is shown that if a connected component $\Gamma$ of the interface is of class $C^2$ and is stationary isothermic, then the mean curvature of $\Gamma$ must be constant. Moreover,  we apply this result to some overdetermined problems for two-phase heat conductors and obtain some symmetry theorems which relax considerably the regularity assumptions of some previous results.
  \end{abstract}



\pagestyle{plain}
\thispagestyle{plain}


\section{Introduction}
\label{introduction}

The paper \cite{KS2021} considered the Cauchy problem for the heat diffusion equation in the whole Euclidean space consisting of two media with different constant conductivities, where initially one medium has temperature 0 and the other  has temperature 1. There, the large time behavior, either stabilization to a constant or oscillation,  of temperature was studied.
The present paper considers the initial behavior of temperature at the interface and studies its interaction with the mean curvature of the interface. 

To be precise, let $\Omega$ be an open set in $\mathbb R^N$ with $N \ge 2$, and let $p \in \partial\Omega$.  Denote by $B_r(x)$ an open ball in $\mathbb R^N$ with radius $r > 0$  and centered at a point $x \in \mathbb R^N$. Assume that there exists  $\rho > 0$ such that $\partial\Omega\cap B_\rho(p)$ is of class $C^2$ and $\partial\Omega$ separates $B_\rho(p)$ into two connected components. Let $\sigma = \sigma(x)\ \left(x \in \mathbb R^N\right)$ be a conductivity distribution of the whole medium satisfying
\begin{equation}
\label{conductivity distribution}
0 < \mu \le \sigma(x) \le M\ \mbox{ for every } x \in \mathbb R^N, \mbox{ and }\sigma(x) = \begin{cases} \sigma_+ \!\! &\mbox{ if }\ x \in B_\rho(p) \cap \Omega,\\  \sigma_-  \!\! &\mbox{ if }\  x \in B_\rho(p) \setminus \Omega,
\end{cases}
\end{equation}
where $\mu, M, \sigma_+ $ and $\sigma_-$ are positive constants. Namely,  the $C^2$ interface $\partial\Omega\cap B_\rho(p)$ separates  $B_\rho(p)$  into two media  with possibly different constant conductivities $\sigma_+, \sigma_-$. In this sense, our setting is purely local in space. Let $u=u(x,t)$ be the unique bounded solution of the Cauchy problem for the heat diffusion equation:
\begin{equation}
\label{Cauchy problem}
u_t=\mbox{div}(\sigma\nabla u)\ \mbox{ in }\ \mathbb R^N\times(0,+\infty)\ \mbox{ and }\ u =\mathcal{X}_{\Omega}\ \mbox{ on }\ \mathbb R^N\times\{0\},
\end{equation}
where $\mathcal X_\Omega$ denotes the characteristic function of the set $\Omega$. The maximum  principle gives
\begin{equation}
\label{positive values}
0 < u(x, t) < 1\ \mbox{ for every } (x,t) \in \mathbb R^N\times (0, + \infty).
\end{equation}
Set
 \begin{equation}
 \label{inside and outside parabolic}
 u^+=u\ \mbox{ if }x \in \overline{\Omega} \quad \mbox{ and } \quad  u^-=u\ \mbox{ if }x \in \mathbb R^N\setminus\Omega.
 \end{equation}
By virtue of \eqref{conductivity distribution},  we have
\begin{equation}
\label{transmission conditions parabolic}
 u^+ =u^-\ \mbox{ and }\ \sigma_+\partial_{\nu}u^+ =  \sigma_-\partial_{\nu}u^-\ \mbox{ on } \left(\partial\Omega\cap B_\rho(p)\right) \times (0, +\infty),
 \end{equation}
where $\nu$ denotes the  outward unit normal vector to $\partial\Omega$ and \eqref{transmission conditions parabolic} expresses the transmission conditions on the $C^2$ interface $\partial\Omega\cap B_\rho(p)$.

It was shown in \cite[Proposition 2.3, p.191]{SaJMPA2020} that, as $t \to 0^+$, the solution $u$ converges to the number $\frac {\sqrt{\sigma_+}}{\sqrt{\sigma_+}+\sqrt{\sigma_-}}$ on $\partial\Omega\cap B_{\rho}(p)$, 
where the convergence is uniform on $\partial\Omega\cap B_r(p)$ for each $0 < r < \rho$. Here the author used the comparison arguments of \cite[Proposition E, p.323]{CMS2021} with the barriers constructed in \cite{Satrieste2016} and the Gaussian bounds for the fundamental solutions of the diffusion equations due to Aronson \cite{Ar1967}. We also prove this convergence in Lemma \ref{new simple barriers} of section \ref{section_a key comparison lemma}  by introducing new simple barriers.  

The first purpose of the present paper is to show the following theorem:

\begin{theorem}
\label{th: asymptotic formula in time}  For every point $x \in \partial\Omega\cap B_{\rho}(p)$, the following formula holds true:
\begin{equation}
\label{asymptotic formula in time}
\lim_{t \to 0^+}t^{-\frac 12}\left(u(x,t)-\frac {\sqrt{\sigma_+}}{\sqrt{\sigma_+}+\sqrt{\sigma_-}}\right) =\frac {\sqrt{\sigma_+}\sqrt{\sigma_-}}{\sqrt{\pi}(\sqrt{\sigma_+}+\sqrt{\sigma_-})} (N-1)H(x),
\end{equation}
where $H(x)$ denotes the mean curvature of $\partial\Omega$ at $x\in\partial\Omega\cap B_\rho(p)$ with respect to the outward normal direction to $\partial\Omega$. The convergence is uniform on $\partial\Omega\cap B_r(p)$ for each $0 < r < \rho$.
\end{theorem}
When $\sigma_+=\sigma_-(=1)$, that is, the equation is just the heat equation, this formula can be directly obtained with the aid of the Gaussian kernel. Indeed, in \cite[Theorem 4.1, pp.546--548]{E1993}, by using  the explicit representation of temperature,  Evans shows  that initially the level surface of temperature $\frac 12$ moves with normal velocity $(N-1)H$.

The second purpose of the present paper  is to give some applications  of Theorem \ref{th: asymptotic formula in time} to some overdetermined problems. The following corollary follows immediately from Theorem \ref{th: asymptotic formula in time}.


\begin{corollary}
\label{constant mean curvature} Let $\sigma=\sigma(x)\ \left(x \in \mathbb R^N\right)$ be given by
\begin{equation}
\label{two-phase}
\sigma(x) = \begin{cases} \sigma_+  \!\!&\mbox{ if }\  x \in  \Omega,\\  \sigma_-  \!\!&\mbox{ if }\  x \in \mathbb R^N \setminus \Omega,
\end{cases}
\end{equation}
and let $u$ be the solution of \eqref{Cauchy problem}. Suppose that $\Gamma$ is a connected component of $\partial\Omega, \ \Gamma$ is of class $C^2$ and there exists a function $a : (0, +\infty) \to (0, +\infty)$ satisfying
\begin{equation}
\label{stationary isothermic}
u(x, t) = a(t)\ \mbox{ for every }\ (x, t) \in \Gamma \times (0, +\infty).
\end{equation}
Then the mean curvature of $\Gamma$ must be constant.
\end{corollary}
In Corollary \ref{constant mean curvature}, the overdetermined condition \eqref{stationary isothermic} means that $\Gamma$ is a stationary isothermic surface.
The following two theorems are examples of the application of Corollary \ref{constant mean curvature}  and Theorem \ref{th: asymptotic formula in time} to the overdetermined problems for two-phase heat conductors where the interface includes a  stationary isothermic surface of class $C^2$.


\begin{theorem}
\label{th:hyperplane}
Let $\Omega \subset \mathbb R^N$ be a domain given by
$$
\Omega = \{ x \in \mathbb R^N :\, x_N > \varphi(x_1,\dots,x_{N-1}) \},
$$
where $\varphi$ is a $C^2$ function on $\mathbb R^{N-1}$. Let $\sigma=\sigma(x)\ \left(x \in \mathbb R^N\right)$ be given by  \eqref{two-phase} and
 let $u$ be the solution of \eqref{Cauchy problem}. Suppose that there exists a function $a : (0, +\infty) \to (0, +\infty)$ satisfying
\begin{equation}
\label{stationary isothermic on entire graph}
u(x, t) = a(t)\ \mbox{ for every }\ (x, t) \in \partial\Omega \times (0, +\infty).
\end{equation}
Then, if $N=2$, $\partial\Omega$ must be a straight line, if $3 \le N \le 8$, $\partial\Omega$ must be a hyperplane, and if $N\ge 9$ and $\nabla \varphi$ is bounded, $\partial\Omega$ must be a hyperplane.
\end{theorem}


\begin{theorem}
\label{th:sphere theorem}
Let $\Omega \subset \mathbb R^N$ be an open set, whose boundary $\partial\Omega$ is of class $C^0$ and has a bounded connected component $\Gamma$ of class $C^2$. Suppose that
either the inside of $\Gamma$ or the outside of $\Gamma$ is included in  one of the two sets, $\Omega$ and $\mathbb R^N\setminus\overline{\Omega}$.
Let $\sigma=\sigma(x)\ \left(x \in \mathbb R^N\right)$ be given by  \eqref{two-phase} and
 let $u$ be the solution of \eqref{Cauchy problem}. Suppose that there exists a function $a : (0, +\infty) \to (0, +\infty)$ satisfying
\begin{equation}
\label{stationary isothermic surface with one C2 bounded part}
u(x, t) = a(t)\ \mbox{ for every }\ (x, t) \in \partial\Omega \times (0, +\infty).
\end{equation}
Then $\partial\Omega$ must be a sphere.
\end{theorem}

Theorem \ref{th:sphere theorem} is also regarded as an improvement of \cite[Theorem 1.1, p.2]{KS2023} where $\partial\Omega$ is assumed to be bounded and of class $C^{2,\alpha}$. Corollary \ref{constant mean curvature} improves the previous two theorems \cite[Theorem 1.1, p.1869]{CSU2023} and \cite[Theorem 1.1, p.187]{SaJMPA2020} in such a way that the assumption that $\partial\Omega$ is uniformly of class $C^6$ can be replaced with that $\partial\Omega$ is uniformly of class $C^2$. Indeed, if $\partial\Omega$ is uniformly of class $C^2$ and stationary isothermic, then Corollary \ref{constant mean curvature} yields that the mean curvature of $\partial\Omega$ is constant and hence $\partial\Omega$ must be real analytic. Furthermore,  the interior estimates for solutions of the prescribed constant mean curvature equation (see \cite[Corollary 16.7, p.407]{GT1983}) give that $\partial\Omega$ is uniformly of class $C^k$ for every $k \in \mathbb N$. Also, by the same reason, in the assumptions of \cite[Theorem 1.5, p.317]{CMS2021}, the $C^6$ regularity of the connected component of $\partial\Omega$ can be replaced with the $C^2$ regularity of it. Namely, Theorem \ref{th: asymptotic formula in time}  and Corollary \ref{constant mean curvature} enable us to relax  considerably the regularity assumptions of some previous results in \cite{KS2023, CSU2023, SaJMPA2020, CMS2021}.

The rest of the paper is organized as follows. 
Section \ref{section_Reduction} is devoted to reducing our setting to the case where $\Omega$ is a bounded $C^2$ domain with two-phase heat conductors as in \eqref{two-phase}. The Gaussian bounds for the fundamental solution of $u_t=\mbox{div}(\sigma\nabla u)$ due to Aronson play a key role.
In section \ref{section_a key comparison lemma}, we prove  a simple  estimate \eqref{new pointwise estimates} in Lemma \ref{new simple barriers} and give the explicit solution \eqref{one dimensional solution}
 in one spatial dimension, both of which play a key role in the proof of Theorem \ref{th: asymptotic formula in time}. 
Section \ref{section_Proof_of_Theorem 1.1} is devoted to the proof of Theorem \ref{th: asymptotic formula in time}. We employ in principle the blow-up arguments due to Ni-Takagi \cite{NT1993} which succeeded to extract the mean curvature of the boundary from the asymptotic behavior of the least-energy solutions of a singularly perturbed semilinear elliptic Neumann problem. One main difference between theirs and ours concerns the limit entire solution appearing in the asymptotics. Namely,  their limit solution is  radially symmetric and  decreasing   but our one is  spatially one-dimensional,  bounded and monotone. Another one is that they consider the homogeneous Neumann condition on the boundary but we consider the transmission conditions on the interface. 

Our proof in section \ref{section_Proof_of_Theorem 1.1} consists of six steps:  We first introduce a principal coordinate system at each point on the interface together with the mean curvature of the interface.  Next we straighten the interface and then introduce the standard parabolic scaling with a small positive parameter $\varepsilon$ for our blow-up arguments. To perform further asymptotics for the scaled solution $v^\varepsilon$ as $\varepsilon \to 0^+$, we utilize the regularity theory, in particular the interior estimates, for the second order parabolic equations of divergence form whose coefficients are H\"older continuous in all but one variable developed by Dong \cite{Dong2012}. Dong's interior estimates guarantee  compactness to employ the blow-up arguments. The most important is the analysis of the second term $S^\varepsilon$ given by \eqref{function for the second mean curvature term} in the asymptotic behavior of $v^\varepsilon$ as $\varepsilon \to 0^+$. Both Lemma \ref{the zeroth-order approximation} coming from Lemma \ref{new simple barriers} and the preliminary estimate \eqref{Lipschitz norm of coefficients for the equation} play a key role in utilizing the interior estimates. Once interior estimates are obtained, the uniqueness of the limit functions guaranteed by the comparison principle for the Cauchy problems plays a key role. Here both of the limit functions $v^*, S^*$ are explicit. Finally, the uniform convergence on the interface in Theorem \ref{th: asymptotic formula in time} also follows from the blow-up arguments.

In section \ref{section_Applications} we prove Theorems  \ref{th:hyperplane} and \ref{th:sphere theorem}. 


\setcounter{equation}{0}
\setcounter{theorem}{0}

\section{Reduction to the case where $\Omega$ is a bounded $C^2$ domain with \eqref{two-phase}}
\label{section_Reduction}

Let $u=u(x,t)$ be the unique bounded solution  of \eqref{Cauchy problem}.  We use the Gaussian bounds for the fundamental solutions of diffusion equations due to
Aronson \cite[Theorem 1, p.891]{Ar1967}(see also \cite[p.328]{FS1986}). Let $g = g(x,\xi,t)$ be the fundamental solution of $u_t=\mbox{div}(\sigma\nabla u)$. Then there exist two positive constants $\kappa< \mathcal K$ depending only on $\mu, M$ and $N$ such that
\begin{equation}
\label{Gaussian bounds}
\kappa\, t^{-\frac N2}e^{-\frac{|x-\xi|^2}{\kappa t}}\le g(x,\xi,t) \le \mathcal K\, t^{-\frac N2}e^{-\frac{|x-\xi|^2}{\mathcal K t}}
\end{equation}
 for all $(x,t), (\xi,t) \in \mathbb R^N\times(0,+\infty)$.  Note that $u$ is represented as
 \begin{equation}
 \label{expression of the solution}
 u(x,t)=\int\limits_\Omega g(x,\xi,t)d\xi\ \mbox{ for } (x,t) \in \mathbb R^N\times(0, +\infty).
 \end{equation}

By the same argument as in \cite[Proof of Proposition 2.2, pp.190--191]{SaJMPA2020} (or as in \cite[Proof of Proposition 2.2, pp.1871--1872]{CSU2023}), we can reduce our setting to the case where $\Omega$ is a bounded $C^2$ domain and  $\sigma$ is given by \eqref{two-phase}. Let $0 < r < \rho$.  Since $\partial\Omega\cap B_\rho(p)$ is of class $C^2$, we may find a bounded domain $\Omega_*$ with $C^2$ boundary $\partial\Omega_*$ satisfying
\begin{equation}
\label{reduction to bounded domain}
\Omega\cap\overline{B_{\frac {r+\rho}2}(p)}\subset \Omega_*\subset \Omega \  \mbox{ and }\  \overline{B_{\frac {r+\rho}2}(p)}\cap\partial\Omega \subset \partial\Omega_*.
\end{equation}
Let us define the conductivity $\sigma_*=\sigma_*(x)\ (x \in \mathbb R^N)$ by
\begin{equation}
\label{two-phase conductivity}
\sigma_*(x) = \begin{cases} \sigma_+ \!\! &\mbox{ if }\ x \in \Omega_*,\\  \sigma_-  \!\! &\mbox{ if }\ x \in \mathbb R^N \setminus \Omega_*.
\end{cases}
\end{equation}
Let $u^*=u^*(x,t)$ be the  unique bounded solution of problem \eqref{Cauchy problem} where $\Omega$ and $\sigma$ are replaced with $\Omega_*$ and $\sigma_*$, respectively. Then, we observe that the difference $w=u-u^*$ satisfies
 \begin{eqnarray}
&&w_t =\mbox{div}(\sigma_*\nabla w)\quad\mbox{in }\ B_{\frac {r+\rho}2}(p)\times (0,+\infty), \label{heat equation initial-boundary*C}
\\ 
&& |w| < 1\  \ \qquad\qquad\mbox{ in }\ \mathbb R^N \times (0,+\infty), \label{heat bounds-2C}
\\
&&w=0  \  \quad\qquad\qquad \mbox{ on } \ B_{\frac {r+\rho}2}(p)\times \{0\}.\label{heat initial*C}
\end{eqnarray}
Set
$$
\mathcal A = \left\{ x \in \mathbb R^N : \mbox{ dist}(x, \partial B_{\frac {r+\rho}2}(p))< \frac 14(\rho-r) \right\} \left(= B_{\frac {3\rho+r}{4}}(p)\setminus \overline{B_{\frac{\rho+3r}{4}}(p)} \right).
$$
By comparing $w$ with the solutions of the Cauchy problem for the heat diffusion equation with conductivity $\sigma_*$ and initial data $\pm 2\mathcal X_{\mathcal A}$ for a short time,  with the aid of the Gaussian bounds \eqref{Gaussian bounds},  we see that
there exist two positive constants $B$ and $b$ such that
\begin{equation}
\label{exponential decay of the difference}
|w(x,t)| \le Be^{-\frac bt}\ \mbox{ for every } (x,t) \in \overline{B_{r}(p) } \times (0,+\infty).
\end{equation}
Therefore,  if  $u^*$ satisfies the conclusion of Theorem \ref{th: asymptotic formula in time},  then $u$ also does.

\setcounter{equation}{0}
\setcounter{theorem}{0}

\section{A key simple estimate and the solution in one dimension}
\label{section_a key comparison lemma}

By virtue of section \ref{section_Reduction}, we may assume that $\Omega$ is a bounded $C^2$ domain and  $\sigma$ is given by \eqref{two-phase}. This assumption is kept through to the end of section \ref{section_Proof_of_Theorem 1.1}.
Let us first prove a new estimate for the solution $u$ of problem \eqref{Cauchy problem} for $x$ in a neighborhood of $\partial\Omega$. This estimate is a key throughout this paper.
We introduce the signed distance function $\delta=\delta(x)$ of $x \in \mathbb R^N$ to the boundary $\partial\Omega$ by
\begin{equation}
\label{signed distance}
\delta(x) = \begin{cases} \mbox{ dist}(x, \partial\Omega)  \!\!&\mbox{ if }\  x \in  \Omega,\\  - \mbox{dist}(x, \partial\Omega)  \!\! &\mbox{ if }\  x \in \mathbb R^N \setminus \Omega.
\end{cases}
\end{equation}
Since $\partial\Omega$ is compact and of class $C^2$, there exists a number $\delta_0 > 0$ such that $\delta(x)$ is of class $C^2$ on the closure of a bounded neighborhood $\mathcal N$ of $\partial\Omega$ given by
\begin{equation}
\label{Neighborhood of the boundary} 
\mathcal N = \{ x \in \mathbb R^N\ :\ -\delta_0 < \delta(x) < \delta_0 \}.
\end{equation}
Define the function $f=f(\xi)$ for $\xi \in \mathbb R$ by
\begin{equation}
\label{error function}
f(\xi)= \frac 1{2\sqrt{\pi}}\int_{-\infty}^\xi e^{-s^2/4}ds.
\end{equation}
Then $f$ satisfies
\begin{eqnarray}
&&f^{\prime\prime}+\frac 12\xi f^\prime=0,\ 0 < f^\prime \le \frac 1{2\sqrt{\pi}} \mbox{ and }\ 0 < f < 1\ \mbox{ in } \mathbb R, \label{equation of selfsimilar solution and monotomicity}\\
&&f(-\infty) =0,\ f(0) = \frac 12,\ f(+\infty) =1,\ f^\prime(0)=\frac 1{2\sqrt{\pi}} \mbox{ and }  f^\prime(\pm\infty)=0.\label{properties of selfsimilar}
\end{eqnarray}
Define the function $\psi=\psi(x, t)\ (x \in \mathbb R^N, t > 0)$ by
\begin{equation}
\label{approximate solutions parabolic}
\psi(x,t) =\begin{cases} 
 \frac {2\sqrt{\sigma_-}}{\sqrt{\sigma_+}+\sqrt{\sigma_-}}\left\{f\!\!\left(\sigma^{-\frac 12}_+t^{-\frac 12}\delta(x)\right) + \frac {\sqrt{\sigma_+}-\sqrt{\sigma_-}}{2\sqrt{\sigma_-}}\right\} & \mbox{if }\ x \in  \Omega,\\
 \frac {2\sqrt{\sigma_+}}{\sqrt{\sigma_+}+\sqrt{\sigma_-}}f\!\!\left(\sigma^{-\frac 12}_-t^{-\frac 12}\delta(x)\right) & \mbox{if }\ x \in  \mathbb R^N\setminus \Omega.
\end{cases}
\end{equation}
Then $\psi$ satisfies the transmission conditions
\begin{equation}
\label{transmission conditions selfsimilar}
\psi|_+=\psi|_- \left(=\frac {\sqrt{\sigma_+}}{\sqrt{\sigma_+}+\sqrt{\sigma_-}}\right)\ \mbox{ and }\ \sigma_+\partial_\nu \psi|_+ =\sigma_-\partial_\nu \psi|_-\ \mbox{ on }\ \partial\Omega \times (0, +\infty),
\end{equation}
where $+$ denotes the limit from the inside of $\Omega$, $-$ does that from the outside of $\Omega$, $\nu =\nu(x)$ denotes the outward unit normal vector to $\partial\Omega$ at $x \in \partial\Omega$ and
$\nu = -\nabla \delta$ on $\partial\Omega$. Since $|\nabla\delta|=1$ in $\mathcal N$, a straightforward computation yields that
\begin{equation}
\label{diffusion equation}
\psi_t-\sigma \Delta \psi =- \frac {2\sqrt{\sigma_+}\sqrt{\sigma_-}}{\sqrt{\sigma_+}+\sqrt{\sigma_-}}t^{-\frac 12}\Delta\delta\  f^\prime\!\!\left(\sigma^{-\frac 12}t^{-\frac 12}\delta(x)\right)\mbox{ in }\  \left(\mathcal N\setminus\partial\Omega\right)\times (0, + \infty).
\end{equation}
We observe from \eqref{Gaussian bounds}, \eqref{expression of the solution}, \eqref{error function} and \eqref{approximate solutions parabolic} that there exists two positive constants $A$ and $a$ such that
\begin{eqnarray}
&0 < u(x,t) \le Ae^{-\frac at}\ &\mbox{ for every } (x,t) \in \left(\partial \mathcal N\setminus\Omega\right)\times (0,+\infty),\label{exponential decay outside Omega for u}\\
&0 < \psi(x,t) \le Ae^{-\frac at}\ &\mbox{ for every } (x,t) \in \left(\partial \mathcal N\setminus\Omega\right)\times (0,+\infty),\label{exponential decay outside Omega for h}\\
&0 < 1-u(x,t) \le Ae^{-\frac at}\ &\mbox{ for every } (x,t) \in \left(\partial \mathcal N\cap\Omega\right)\times (0,+\infty),\label{exponential decay inside Omega for u}\\
&0 < 1-\psi(x,t) \le Ae^{-\frac at}\ &\mbox{ for every } (x,t) \in \left(\partial \mathcal N\cap\Omega\right)\times (0,+\infty).\label{exponential decay inside Omega for h}
\end{eqnarray}

\begin{lemma}
\label{new simple barriers} Let $u$ be the solution of problem \eqref{Cauchy problem}.
Then, there exists a positive constant $K$ such that
\begin{equation}
\label{new pointwise estimates}
\psi(x,t)-K\sqrt{t} \le u(x,t) \le \psi(x,t)+K\sqrt{t}\ \mbox{ for every } (x,t) \in \mathcal N \times (0,+\infty).
\end{equation}
Moreover, as $t \to 0^+$, $u$ converges to the number  $\frac {\sqrt{\sigma_+}}{\sqrt{\sigma_+}+\sqrt{\sigma_-}}$ uniformly on $\partial\Omega$.
\end{lemma}

\noindent
{\it Proof.\ } Since $\psi=\frac {\sqrt{\sigma_+}}{\sqrt{\sigma_+}+\sqrt{\sigma_-}}$ on $\partial\Omega$, the latter follows immediately from \eqref{new pointwise estimates}. Let us prove the former.
Set $K=\max\{ \frac {4\sqrt{\sigma_+}\sqrt{\sigma_-}}{\sqrt{\sigma_+}+\sqrt{\sigma_-}}\max\limits_{\overline{\mathcal N}}|\Delta\delta|  \max\limits_{\xi\in \mathbb R}f^\prime(\xi), \frac A{\sqrt{2ea}}\}$ and define the two functions $k_\pm=k_\pm(x,t)$ by
$$
k_\pm(x,t) = \psi(x,t)\pm K\sqrt{t}\ \mbox{ for } (x,t) \in \overline{\mathcal N} \times (0, + \infty).
$$
Then it follows from \eqref{transmission conditions selfsimilar} and \eqref{diffusion equation} that $k_\pm$ satisfy the following in the weak sense:
\begin{equation}
\label{sub and super solutions}(k_-)_t - \mbox{div}(\sigma\nabla k_-) \le 0 \le (k_+)_t - \mbox{div}(\sigma\nabla k_+) \ \mbox{ in }\ \mathcal N \times (0, +\infty).
\end{equation}
Moreover, we observe from \eqref{exponential decay outside Omega for u}--\eqref{exponential decay inside Omega for h} that 
\begin{equation}
\label{boundary condition}
k_-\le u\le k_+\ \mbox{ on } \partial\mathcal N \times (0, +\infty).
\end{equation}
Therefore, since $u = k_+ =k_- =\mathcal X_\Omega$ on $\mathcal N \times \{0\}$, we have from \eqref{sub and super solutions},  \eqref{boundary condition} and the comparison principle that
$$
k_-\le u\le k_+\  \mbox{ in }\mathcal N \times (0, +\infty),
$$
which is the conclusion. \qed

Finally in this section, we mention some useful properties of the solution in one spatial dimension.
Define the function $\Psi=\Psi(\eta, s)\ (\eta \in \mathbb R, s > 0)$ by
\begin{equation}
\label{one dimensional solution}
\Psi(\eta, s) =\begin{cases} 
 \frac {2\sqrt{\sigma_-}}{\sqrt{\sigma_+}+\sqrt{\sigma_-}}\left\{f\!\!\left(\sigma^{-\frac 12}_+s^{-\frac12}\eta\right) + \frac {\sqrt{\sigma_+}-\sqrt{\sigma_-}}{2\sqrt{\sigma_-}}\right\} & \mbox{if }\ \eta > 0,\\
 \frac {2\sqrt{\sigma_+}}{\sqrt{\sigma_+}+\sqrt{\sigma_-}}f\!\!\left(\sigma^{-\frac 12}_- s^{-\frac 12}\eta\right) & \mbox{if }\ \eta \le 0,
\end{cases}
\end{equation}
 where $f$ is given by \eqref{error function}.
Note that 
\begin{equation}
\label{Psi and psi relationship}
\psi(x,t)=\Psi(\delta(x), t)\ \mbox{ for } (x,t)\in\mathbb R^N\times(0,+\infty).
\end{equation}
Then it follows from the definition  \eqref{error function} of $f$  that the transmission conditions at $(0, s)$ for every $s > 0$ hold true:
\begin{equation}
\label{one-dimensional transmission}
\Psi(+0, s)=\Psi(-0, s)\left(=  \frac {\sqrt{\sigma_+}}{\sqrt{\sigma_+}+\sqrt{\sigma_-}}\right)\ \mbox{ and }\  \sigma_+\partial_\eta\Psi(+0,s)=\sigma_-\partial_\eta\Psi(-0,s).
\end{equation}
Moreover, we have from  \eqref{equation of selfsimilar solution and monotomicity} and \eqref{properties of selfsimilar} that
\begin{equation}
\label{properties of Phi}
0 < \Psi < 1 \mbox{ and } \partial_\eta\Psi > 0\ \mbox{ in } \mathbb R\times (0,+\infty),\ \Psi(-\infty,s) = 0,\ \Psi(+\infty, s)=1\ \mbox{ for } s > 0.
\end{equation}
A straightforward computation with the aid of  \eqref{equation of selfsimilar solution and monotomicity} gives
\begin{equation}
\label{equations selfsimilar one-dimensional}
\Psi_s= \begin{cases} \sigma_+\partial_\eta^2\,\Psi\ \!\!&\mbox{ if }\  \eta >0,\\  \sigma_- \partial_\eta^2\,\Psi\ \!\!&\mbox{ if }\  \eta < 0. \end{cases}
\end{equation}
 Let us introduce two functions $\sigma^*$ and $\mathcal X^*$ by
\begin{equation}
\label{one dimensional conductivity distribution and characteristic function}
\sigma^*=\sigma^*(\eta) = \begin{cases} \sigma_+ \!\!&\mbox{ if }\  \eta >0,\\  \sigma_-  \!\!&\mbox{ if }\  \eta \le 0, \end{cases}\ \mbox{ and }\ \mathcal X^*=\mathcal X^*(\eta)= \begin{cases} 1  \!\!&\mbox{ if }\  \eta >0,\\  0  \!\!&\mbox{ if }\  \eta \le 0.
\end{cases}
\end{equation}
Then, by combining \eqref{one-dimensional transmission}--\eqref{one dimensional conductivity distribution and characteristic function} with the fact that $f(-\infty) =0$ and $f(+\infty)=1$ in \eqref{properties of selfsimilar},  we see that  $\Psi$ is the unique bounded weak solution of the Cauchy problem:
\begin{equation}
\label{one dimensional entire bounded solution of Cauchy}
\Psi_s =\partial_\eta\!\left(\sigma^* \partial_\eta\Psi\right) \mbox{ in } \mathbb R \times (0,+\infty) \ \mbox{ and }\ \Psi = \mathcal X^* \mbox{ on } \mathbb R\times\{0\}.
\end{equation}

\setcounter{equation}{0}
\setcounter{theorem}{0}

\section{Proof of Theorem \ref{th: asymptotic formula in time}}
\label{section_Proof_of_Theorem 1.1}

We employ in principle  the blow-up arguments due to Ni-Takagi \cite{NT1993} which succeeded to extract the mean curvature of the boundary from the asymptotic behavior of the least-energy solutions of a singularly perturbed semilinear  elliptic Neumann problem.
By straightening the boundary and introducing the scaling related to the small parameter $\varepsilon$ in their singularly perturbed problem, they find the mean curvature of the boundary in the elliptic equation satisfied by the first-order approximation in $\varepsilon$ as $\varepsilon \to 0^+$ (see \cite[(2.4), p.251 and (2.8), p.252]{NT1993}). 

Here we straighten the interface and introduce the standard parabolic scaling with a small positive parameter $\varepsilon$. Then we can find the mean curvature of the interface in the inhomogeneous diffusion equation \eqref{limit equation 2nd} satisfied by the first-order approximation $S^*$ in $\varepsilon$  as $\varepsilon \to 0^+$ of the scaled solution $v^\varepsilon$. Our proof consists of six steps.

\subsection{Introducing a principal coordinate system at each point on the interface}
\label{subsection principal coordinate system}
Let $q \in \partial\Omega$ arbitrarily. For each point $q$, let us introduce a principal coordinate system with the origin $0$ at $q\in\partial\Omega$. Since $\partial\Omega$ is compact and of class $C^2$, we may choose a number $r_0$ with $0 < r_0< \delta_0$, which is independent of $q$, such that there exists a function $\varphi \in C^2(\mathbb R^{N-1})$ satisfying
$$
B_{r_0}(q)\cap\Omega=\{  x \in B_{r_0}(0) : x_N < \varphi(\hat{x}) \} \mbox{ and } B_{r_0}(q)\cap\partial\Omega=\{ x \in B_{r_0}(0) : x_N=\varphi(\hat{x})\},
$$
where $\delta_0 $ is the positive constant given in \eqref{Neighborhood of the boundary}, $\hat{x}=(x_1,\dots,x_{N-1})$ for $x \in \mathbb R^N$, and the norm $\Vert\varphi\Vert_{C^2(\mathbb R^{N-1})}$  is bounded in $q \in \partial\Omega$.
Without loss of generality, we may suppose that 
\begin{equation}
\label{principal curvatures}
q=0,  \varphi(0)=0, \nabla\varphi(0)=0\mbox{ and }\varphi(\hat{x})=\frac 12\sum_{j=1}^{N-1}\kappa_j x_j^2 + o(|\hat{x}|^2)\ \mbox{ as } \hat{x} \to 0,
\end{equation}
where $\kappa_j\ (j=1,\dots,N\!-\!1)$ are the principal curvatures of $\partial\Omega$ at $q \in \partial\Omega$ with respect to the outward normal direction to $\partial\Omega$.  Note that the mean curvature $H(q)$ of
 $\partial\Omega$ at $q\in\partial\Omega$ is given by
 $$
 H(q)=H(0)=\frac 1{N-1}\sum\limits_{j=1}^{N-1}\kappa_j.
 $$ 
 
 
 \subsection{Straightening the interface and scaling for blow-up arguments}
 \label{subsection straightening the interface and scaling}
 As used in the boundary regularity theory for elliptic partial differential equations in \cite{ACM2018, E2010}, let us introduce the coordinate transformation $T : (\hat{x}, x_N) \mapsto (\hat{y}, y_N)=(\hat{x}, \varphi(\hat{x})-x_N)$ which straightens $B_{r_0}(q)\cap\partial\Omega$.  Then,  the Jacobian matrix $\nabla T$ is triangular with $\mbox{det}(\nabla T)=-1$,  its inverse transformation is given by $T^{-1}:  (\hat{y}, y_N) \mapsto (\hat{x}, x_N)=(\hat{y},\varphi(\hat{y})-y_N)$, and $\nabla T^{-1}$ is also triangular with $\mbox{det}(\nabla T^{-1})=-1$. We may put 
 $$
 D=T(B_{r_0}(0)),\ T(B_{r_0}(0)\cap\Omega)=\{ y \in D :  y_N > 0\} \mbox{ and }T(B_{r_0}(0)\cap\partial\Omega)= \{ y \in D : y_N = 0\}.
 $$
 Then,  by setting for $(y,t) \in D\times [0, +\infty)$ 
$$
 v=v(y,t) = u(T^{-1}y,t),\ \sigma_T=\sigma_T(y_N) = \sigma(T^{-1}y)\ \mbox{ and } \mathcal X_T= \mathcal X_T(y_N) = \mathcal X_\Omega(T^{-1}y),
 $$
we see from \eqref{Cauchy problem} and \eqref{positive values} that
\begin{eqnarray}
 &0 < v < 1\ \mbox{ and }\ v_t=\mbox{div}(\sigma_T A\nabla v)\ &\mbox{ in } D \times (0,+\infty), \label{weak solution after straightening}\\
 &v= \mathcal X_{T} &\mbox{ on } D \times \{0\}, \label{initial data after straightening}
 \end{eqnarray}
where  $A=A(\hat{y})$ is the symmetric $N\times N$ matrix given by
\begin{equation}
\label{symmetric uniformly elliptic matrix}
A(\hat{y}) = \left[\begin{array}{c|c} I_{N-1}&\nabla \varphi(\hat{y})\\ \hline
{}^t(\nabla\varphi(\hat{y}))&1+|\nabla\varphi(\hat{y})|^2
\end{array}\right]\ \mbox{ and }\ A(0)=I_N
\end{equation}
 with the identity matrix $I_{k}$ of size $k \in \mathbb N$ and the transpose ${}^t(\nabla\varphi(\hat{y}))$ of  $\nabla \varphi(\hat{y})$.

 We observe that 
 $$
 \sigma_T(y_N) = \sigma^*(y_N)\ \mbox{ and }\mathcal X_{T}(y_N) = \mathcal X^*(y_N)\ \mbox{  for }y \in D,
 $$
 where $\sigma^*$ and $\mathcal X^*$ are given by \eqref{one dimensional conductivity distribution and characteristic function}.

Since the norm $\Vert\varphi\Vert_{C^2(\mathbb R^{N-1})}$  is bounded in $q \in \partial\Omega$ and $T0=0$, there exists a number $\rho_0 > 0$ being independent of $q \in \partial\Omega$ such that $\overline{B_{\rho_0}(0)} \subset D$. Let us  introduce the standard parabolic scaling with a small positive parameter $\varepsilon$ by 
$$
(z,s)=(\varepsilon^{-1}\,y, \varepsilon^{-2}t)\ \mbox{ for }(y, t) \in \overline{B_{\rho_0}(0)} \times [0,+\infty).
$$
Remark that $y \in B_{\rho_0}(0)$ if and only if $z\in B_{\varepsilon^{-1}\rho_0}(0)$. For $(z,s) \in B_{\varepsilon^{-1}\rho_0}(0)\times [0, +\infty)$, we set
\begin{eqnarray}
&&v^\varepsilon=v^\varepsilon(z,s) = v\!\left(\varepsilon z, \varepsilon^2 s \right),\  A^*=A^*(\hat{z}) = A\!\left(\varepsilon\hat{z}\right), \label{scaled solution}\\
&&\sigma^*= \sigma^*(z_N)\ \mbox{ and }\ 
\mathcal X^*=\mathcal X^*(z_N),\label{scaled conductivities etc}
\end{eqnarray}
where $\sigma^*, \mathcal X^*$ are given by \eqref{one dimensional conductivity distribution and characteristic function} and these  are invariant under the scaling.
Then it follows from \eqref{weak solution after straightening} and \eqref{initial data after straightening} that
\begin{eqnarray}
 &0 < v^\varepsilon < 1\ \mbox{ and }\ v^\varepsilon_s=\mbox{div}(\sigma^* A^*\nabla v^\varepsilon)\ &\mbox{ in } B_{\varepsilon^{-1}\rho_0}(0)  \times (0,+\infty), \label{weak solution after scaling}\\
 &v^\varepsilon= \mathcal X^*&\mbox{ on }  B_{\varepsilon^{-1}\rho_0}(0)\times \{0\}. \label{initial data after scaling}
 \end{eqnarray}

\subsection{Utilizing the interior estimates for two-phase problems}
\label{subsection interior estimates}
Let us  utilize the regularity theory for the second order parabolic equations of divergence form whose coefficients are H\"older continuous in all but one variable developed by \cite{Dong2012}. See \cite{CKV1986, Dong2012, LN2003, LV2000, XB2013} for elliptic equations.  In our setting \eqref{weak solution after scaling}--\eqref{initial data after scaling},  this one variable means $z_N$ and  Lipschitz continuity in $\hat{z}$ of the coefficients $\sigma^* A^*$ is guaranteed by \eqref{symmetric uniformly elliptic matrix} with $\varphi \in C^2(\mathbb R^{N-1})$. The discontinuity in $z_N$ of  $\sigma^* A^*$ comes from that of $\sigma^*$ at $z_N=0$. Note that if $n \in \mathbb N$ and  $\varepsilon \le (n+2)^{-1}$, then $\overline{B_{n\rho_0}(0)} \subset B_{(n+1)\rho_0}(0)\subset\overline{B_{(n+1)\rho_0}(0)}\subset B_{\varepsilon^{-1}\rho_0}(0)$. For each $n \in \mathbb N$, we define a cylinder $Q_n$ in $\mathbb R^{N}\!\times\!(0,+\infty)$ by
\begin{equation}
\label{sequence of cylinders}
Q_n= B_{n\rho_0}(0)\times(n^{-1}\rho_0^2, n^{-1}\rho_0^2+n^2\rho_0^2).
\end{equation}
Hence $\overline{Q_{n}}  \subset  Q_{n+1}\subset\overline{Q_{n+1}}\subset B_{\varepsilon^{-1}\rho_0}(0)\times (0, +\infty)$ and $\bigcup\limits_{n=1}^\infty Q_n = \mathbb R^N\times(0,+\infty)$.
 Let $0 < \alpha < 1$. Then,  with the aid of Dong's interior estimates \cite[Theorem 2, p.125]{Dong2012}, we infer that, for every $n \in \mathbb N$ there exists a positive constant $C_n$, which is independent of $q \in \partial\Omega$, such that  for every $0 <\varepsilon \le (n+2)^{-1}$
\begin{equation}
\label{interior estimates}
 \Vert v^\varepsilon\Vert_{C^{\alpha,\alpha/2}(\overline{Q_n})} + \Vert\nabla_{\hat{z}} v^\varepsilon\Vert_{C^{\alpha,\alpha/2}(\overline{Q_n})}+ \Vert \partial_{z_N}\! v^\varepsilon\Vert_{C^{\alpha,\alpha/2}(\overline{Q_{n+}})}+ \Vert \partial_{z_N}\! v^\varepsilon\Vert_{C^{\alpha,\alpha/2}(\overline{Q_{n-}})} \le C_n,
 \end{equation}
 where $Q_{n\pm}= \{ (z,s) \in Q_n : \pm z_N > 0 \},\ \nabla_{\hat{z}} v^\varepsilon$ denotes the gradient in $\hat{z} \in \mathbb R^{N-1}$, and the norm $ \Vert w \Vert_{C^{\alpha,\alpha/2}(E)}$ for a function $w=w(z,s)$ on $E \subset \mathbb R^{N+1}$  is given by 
 \begin{equation}
 \label{Holder norm definition}
 \Vert w \Vert_{C^{\alpha,\alpha/2}(E)} = \sup_E |w| +\sup_{\substack{(z,s), (\tilde{z},\tilde{s}) \in E\\ (z,s)\not= (\tilde{z},\tilde{s})}}\frac {|w(z,s)-w(\tilde{z},\tilde{s})|}{|z-\tilde{z}|^\alpha+|s-\tilde{s}|^{\alpha/2}}.
 \end{equation}
 The interior estimates obtained by \cite[Theorem 2, p.125]{Dong2012} are sharper than \eqref{interior estimates}, but for our purpose \eqref{interior estimates} is sufficient.

 \subsection{Getting the first term in the asymptotic formula}
 \label{subsection the first term}
  By virtue of \eqref{interior estimates}, it follows from the Arzel\`a-Ascoli theorem together with the Cantor diagonal process that there exist a sequence $\{\varepsilon_j\}$ with $\lim\limits_{j\to \infty}\varepsilon_j=0$ and a function $v^*$ on $\mathbb R^N\times(0,+\infty)$ which satisfy the following for every $n \in \mathbb N$:
  \begin{eqnarray}
  \{v^{\varepsilon_j}\}_{j\ge n} \mbox{ converges to } v^* \mbox{ as } j \to \infty \mbox{ uniformly on } \overline{Q_n};\ \ &&\label{uniform convergence subsequence 1st}\\
   \{\nabla v^{\varepsilon_j}\}_{j\ge n} \mbox{ converges to } \nabla v^* \mbox{ as } j \to \infty \mbox{ uniformly on } \overline{Q_{n\pm}};\ \ &&\label{uniform convergence subsequence for gradient 1st}
  \\
  \Vert v^*\!\Vert_{C^{\alpha,\alpha/2}(\overline{Q_n})} + \Vert\nabla_{\hat{z}} v^*\!\Vert_{C^{\alpha,\alpha/2}(\overline{Q_n})}+ \Vert \partial_{z_N}\! v^*\!\Vert_{C^{\alpha,\alpha/2}(\overline{Q_{n+}})}+ \Vert \partial_{z_N}\! v^*\!\Vert_{C^{\alpha,\alpha/2}(\overline{Q_{n-}})} \le C_n,\ \ &&\label{estimate for limit function 1st}.
  \end{eqnarray}
   where  $C_n, \alpha$ are the constants given in \eqref{interior estimates}. Moreover, by \eqref{weak solution after scaling}, $v^*$ satisfies
  \begin{equation}
  0 \le  v^*\le 1\ \mbox{ and } v_s^*=\mbox{div}(\sigma^* \nabla v^*)  \ \mbox{ in } \mathbb R^N\times (0,+\infty).\label{limit equation 1st} 
 \end{equation}
 Here  we used the fact that, as $\varepsilon\to 0^+$, $A^*$ converges to $A(0) = I_N$ uniformly on every compact set in $\mathbb R^N$. 
 
 Next lemma shows that $v^*(z,s)$ is explicitly given by $\Psi(z_N,s)$ and it must be the zeroth-order approximation in $\varepsilon$ as $\varepsilon \to 0^+$ of $v^\varepsilon$.
 \begin{lemma}
 \label{the zeroth-order approximation}
 The following assertions hold:
 \begin{itemize}
 \item[\rm (i)] For each $n\in\mathbb N$, as $\varepsilon \to 0^+$, $|v^\varepsilon(z,s)- \Psi(z_N,s)|\le \varepsilon(K\sqrt{s} +o(1))$ for $(z,s)\in Q_n$;
 \item[\rm (ii)] $v^*(z,s) = \Psi(z_N,s)$ for all $(z,s) \in \mathbb R^N\times(0,+\infty)$;
 \item[\rm (iii)] For each $n\in\mathbb N$, $ \{v^{\varepsilon}\}_{\varepsilon\le (n+2)^{-1}}$ converges to $ v^*$ as  $\varepsilon \to 0^+$ uniformly on $\overline{Q_n}$;
  \item[\rm (iv)]For each $n\in\mathbb N$,  $ \{ \nabla v^{\varepsilon}\}_{\varepsilon\le (n+2)^{-1}}$  converges to $\nabla v^*$ as $ \varepsilon \to 0^+$ uniformly on $\overline{Q_{n\pm}}$,
 \end{itemize}
 where $K$ is the constant given by Lemma \ref{new simple barriers} and $\Psi$ is the function given by \eqref{one dimensional solution}.
 \end{lemma}
 
 \noindent
 {\it Proof. } (i) together with \eqref{uniform convergence subsequence 1st} yields (ii). Both (iii) ad (iv) follow from (ii) with the aid of the compactness arguments in the beginning of  subsection \ref{subsection the first term}. It remains to prove (i).
  Let $n \in \mathbb N$ and $0 < \varepsilon \le (n+2)^{-1}$. Then $\overline{B_{n\rho_0}(0)}\subset B_{\varepsilon^{-1}\rho_0}(0)$. The signed distance function $\delta=\delta(x)$ given by \eqref{signed distance} satisfies
\begin{equation}
\label{Taylor expansions of the signed distance function}
\delta(x) = -x_N + \frac 12\sum_{j=1}^{N-1}\kappa_j x_j^2 + o(|x|^2) \mbox{ as } x \to 0.
\end{equation}
Since $\hat{x}=\varepsilon\hat{z}$ and $x_N=\varphi\!\left(\varepsilon\hat{z}\right)-\varepsilon z_N$, we observe that, for each $z \in B_{n\rho_0}(0)$,  
\begin{equation}
\label{asymptotics of xN}
x_N=\frac 12\sum_{j=1}^{N-1}\kappa_j\!\left(\varepsilon z_j\right)^2-\varepsilon z_N+o\!\left(\varepsilon^2\right)\ \mbox{ as } \varepsilon \to 0^+.
\end{equation}
Hence,  $|x|^2=\varepsilon^2|z|^2 + O\!\left(\varepsilon^3\right)  \mbox{ as } \varepsilon \to 0^+$, and combining  \eqref{Taylor expansions of the signed distance function} with \eqref{asymptotics of xN} yields that for each $z \in B_{n\rho_0}(0)$
\begin{equation}
\label{a key estimate}
\delta(x) = \varepsilon z_N + o\!\left(\varepsilon^2\right)\ \mbox{ as } \varepsilon \to 0^+.
\end{equation}
Let $(z,s)\in Q_n$. Since $t=\varepsilon^2s$, it follows from \eqref{a key estimate} that
\begin{equation}
\label{distance function with time}
t^{-\frac 12}\delta(x)=\varepsilon^{-1} s^{-\frac 12}(\varepsilon z_N + o\!\left(\varepsilon^2\right)) =s^{-\frac 12}z_N+o(\varepsilon)\ \mbox{ as } \varepsilon \to 0^+.
\end{equation}
Then, by recalling  that the function $f$ given by \eqref{error function} satisfies $0< f^\prime\le\frac 1{2\sqrt{\pi}}$, we have
\begin{equation}
\label{difference of f}
f(\sigma_\pm^{-\frac 12}t^{-\frac 12}\delta(x))-f(\sigma_\pm^{-\frac 12}s^{-\frac12}z_N)=o(\varepsilon)\ \mbox{ as }\varepsilon \to 0^+,
\end{equation}
 which implies that
 \begin{equation}
\label{difference of psi  and Psi}
\psi(x,t)-\Psi(z_N,s)=o(\varepsilon)\ \mbox{ as }\varepsilon \to 0^+,
\end{equation}
where $\psi, \Psi$ are given by \eqref{approximate solutions parabolic}, \eqref{one dimensional solution}, respectively. Therefore, since $u(x,t)=v^\varepsilon(z,s)$, we have from Lemma \ref{new simple barriers} that as $\varepsilon \to 0^+$
$$
|v^\varepsilon(z,s)- \Psi(z_N,s)|\le|u(x,t)-\psi(x,t)|+|\psi(x,t)- \Psi(z_N,s)| \le K\sqrt{t} + o(\varepsilon) =\varepsilon K\sqrt{s} +o(\varepsilon),
$$
 which gives (i). \qed


 \subsection{Getting  the second term in the asymptotic formula}
 \label{subsection the second term}
Let us introduce the function $S^\varepsilon=S^\varepsilon(z, s)$ for $(z,s) \in B_{\varepsilon^{-1}\rho_0}(0)  \times (0,+\infty)$ by
\begin{equation}
\label{function for the second mean curvature term}
S^\varepsilon(z,s) = \varepsilon^{-1}(v^\varepsilon(z,s)-v^*(z,s)) \ \left\{ = \varepsilon^{-1}(v^\varepsilon(z,s)-\Psi(z_N,s)) \right\}.
\end{equation}
Then we have from (i) of Lemma \ref{the zeroth-order approximation} that for each $n\in\mathbb N$, as $\varepsilon \to 0^+$, 
\begin{equation}
\label{bounds of the second term}
|S^\varepsilon(z, s)|\le K\sqrt{s} +o(1)\ \mbox{ for } (z,s)\in Q_n.
\end{equation}

Let us derive the equation which $S^\varepsilon$ satisfies and the asymptotic behavior of $S^\varepsilon$ as $\varepsilon \to 0^+$. Since it follows from \eqref{function for the second mean curvature term} that
$$
v^\varepsilon(z,s) = \varepsilon S^\varepsilon(z,s) + \Psi(z_N,s),
$$
we see from \eqref{weak solution after scaling}  and \eqref{one dimensional entire bounded solution of Cauchy} that 
\begin{equation}
\label{weak solution to the second order approximate}
S^\varepsilon_s=\mbox{div}(\sigma^*A^*\nabla S^\varepsilon) + \mbox{div}(\varepsilon^{-1}\!\sigma^*(A^*-I_N)\nabla \Psi)\ \mbox{ in } B_{\varepsilon^{-1}\rho_0}(0)\times (0,+\infty).
 \end{equation}
Let us compute the second term of the right-hand side.  We observe that for $(z,s) \in B_{\varepsilon^{-1}\rho_0}(0)\times (0,+\infty)$ 
\begin{eqnarray*}
\sigma^*\nabla \Psi(z_N,s) &=& {}^t[0,\dots,0,\sigma^*\partial_{z_N}\Psi(z_N,s)],\\
A^*(z)-I_N&=&A\!\left(\varepsilon\hat{z}\right)-I_N= \left[\begin{array}{c|c} O&(\nabla \varphi)(\varepsilon\hat{z})\\ \hline
{}^t((\nabla\varphi)(\varepsilon\hat{z}))&|(\nabla\varphi)(\varepsilon\hat{z})|^2
\end{array}\right].
\end{eqnarray*}
Let $n \in \mathbb N, \ \varepsilon \le (n+2)^{-1}$ and $(z,s) \in Q_{n+1}$. Remark that $\overline{B_{(n+1)\rho_0}(0)}\subset B_{\varepsilon^{-1}\rho_0}(0)$.
Since  $\Vert\varphi\Vert_{C^2(\mathbb R^{N-1})}$ is finite and $\nabla\varphi(0)=0$,   for every $j=1,\dots,N\!-\!1$ and every $\varepsilon \le (n+2)^{-1}$, we have the following estimates:
\begin{eqnarray*}
&&\sup_{z \in \overline{B_{(n+1)\rho_0}(0)}}\!\!\varepsilon^{-1}|(\nabla\varphi)(\varepsilon\hat{z})| \le (n+1)\rho_0\, \Vert\varphi\Vert_{C^2(\mathbb R^{N-1})}\\
 &&\mbox{ and } \sup_{z \in \overline{B_{(n+1)\rho_0}(0)}}\!\!\varepsilon^{-1}|\partial_{z_j}\{(\nabla\varphi)(\varepsilon\hat{z})\}|\le \Vert\varphi\Vert_{C^2(\mathbb R^{N-1})}.
\end{eqnarray*}
Thus, since $\sigma^*\partial_{z_N}\Psi \in C^{1,1/2}(\overline{Q_{n+1}})$ because of \eqref{one dimensional solution}, we infer that there exists a constant $C_{1,n} >0$, which is independent of $q \in \partial\Omega$, satisfying for every $\varepsilon\le(n+2)^{-1}$
\begin{equation}
\label{Lipschitz norm of coefficients for the equation}
\Vert \varepsilon^{-1}\, \sigma^*(A^*-I_N)\nabla \Psi\Vert_{C^{1,1/2}(\overline{Q_{n+1}})} \le C_{1,n},
\end{equation}
where the norm of the left-hand side means \eqref{Holder norm definition} with $\alpha=1$.
By virtue of \eqref{bounds of the second term} for $n+1$, \eqref{weak solution to the second order approximate} and \eqref{Lipschitz norm of coefficients for the equation}, we can again use  Dong's interior estimates \cite[Theorem 2, p.125]{Dong2012}.
Then there exists a constant $C_{2,n} > 0$ being independent of $q \in \partial\Omega$ such that for every $\varepsilon\le(n+2)^{-1}$
\begin{equation}
\label{interior estimates for S varepsilon}
 \Vert S^\varepsilon\Vert_{C^{\alpha,\alpha/2}(\overline{Q_n})} + \Vert\nabla_{\hat{z}} S^\varepsilon\Vert_{C^{\alpha,\alpha/2}(\overline{Q_n})}+ \Vert \partial_{z_N}\! S^\varepsilon\Vert_{C^{\alpha,\alpha/2}(\overline{Q_{n+}})}+ \Vert \partial_{z_N}\! S^\varepsilon\Vert_{C^{\alpha,\alpha/2}(\overline{Q_{n-}})} \le C_{2,n},\qquad
 \end{equation}
 where  $\alpha\in (0,1)$ is the same number as in \eqref{interior estimates} which is independent of $n\in\mathbb N$ and $q \in \partial\Omega$.  Therefore, by virtue of \eqref{interior estimates for S varepsilon}, it follows from the Arzel\`a-Ascoli theorem together with the Cantor diagonal process that there exist a sequence $\{\varepsilon_j\}$ with $\lim\limits_{j\to \infty}\varepsilon_j=0$ and a function $S^*$ on $\mathbb R^N\times (0,+\infty)$ which satisfy the following for every $n \in \mathbb N$:
 \begin{eqnarray*}
  \{S^{\varepsilon_j}\}_{j\ge n} \mbox{ converges to } S^* \mbox{ as } j \to \infty \mbox{ uniformly on } \overline{Q_n}; &&\\
   \{\nabla S^{\varepsilon_j}\}_{j\ge n} \mbox{ converges to } \nabla S^* \mbox{ as } j \to \infty \mbox{ uniformly on } \overline{Q_{n\pm}}; &&
  \\
  \Vert S^*\!\Vert_{C^{\alpha,\alpha/2}(\overline{Q_n})} + \Vert\nabla_{\hat{z}} S^*\!\Vert_{C^{\alpha,\alpha/2}(\overline{Q_n})}+ \Vert \partial_{z_N}\! S^*\!\Vert_{C^{\alpha,\alpha/2}(\overline{Q_{n+}})}+ \Vert \partial_{z_N}\! S^*\!\Vert_{C^{\alpha,\alpha/2}(\overline{Q_{n-}})} \le C_{2,n}, &&
  \end{eqnarray*}
   where  $C_{2,n}, \alpha$ are the same constants as in \eqref{interior estimates for S varepsilon}.   Moreover, since for $i=1,\dots,N\!-\!1$ 
$$
(\partial_{x_i} \varphi)(\varepsilon\hat{z}) = \sum_{k=1}^{N-1}(\partial_{x_k} \partial_{x_i} \varphi)(0)\varepsilon z_k + o(\varepsilon)
=\kappa_i\varepsilon z_i + o(\varepsilon)\ \mbox{ as } \varepsilon \to 0^+,
$$
$\varepsilon^{-1}\, \sigma^*(A^*-I_N)\nabla \Psi$ converges to $\sigma^*\partial_{z_N}\Psi\,{}^t[\kappa_1z_1,\dots,\kappa_{N-1}z_{N-1}, 0]$ as $\varepsilon \to 0^+$ uniformly  on every compact set in $\mathbb R^N\times (0,+\infty)$. 
Thus we infer from \eqref{bounds of the second term} and \eqref{weak solution to the second order approximate} that
  \begin{eqnarray}
 && |S^*(z,s)| \le K\sqrt{s}\ \mbox{ for every } (z,s) \in \mathbb R^N\times(0,+\infty); \\
 && S^*_s =\mbox{div}(\sigma^* \nabla S^*) +(N-1)H(0)\sigma^*\partial_{z_N}\Psi \ \mbox{ in } \mathbb R^N\times(0, +\infty), \label{limit equation 2nd} 
 \end{eqnarray}
 where  we also used the fact that, as $\varepsilon\to 0^+$, $A^*$ converges to $A(0) = I_N$ uniformly on every compact set in $\mathbb R^N$. These guarantee that, for each $\tau > 0$, $S^*$ must be the unique bounded solution of the Cauchy problem:
 $$
 S^*_s =\mbox{div}(\sigma^* \nabla S^*) +(N-1)H(0)\sigma^*\partial_{z_N}\Psi \ \mbox{ in } \mathbb R^N\times(0, \tau] \ \mbox{ and }\ S^*=0\ \mbox{ on }  \mathbb R^N\times\{0\}.
 $$
 Furthermore, this is explicitly solved and hence we have
 \begin{equation}
 \label{the explicit solution to the second approximate}
 S^*=S^*(z_N,s)=\frac {2\sqrt{\sigma_+}\sqrt{\sigma_-}}{\sqrt{\sigma_+}+\sqrt{\sigma_-}}(N-1)H(0)\sqrt{s}\begin{cases} f^\prime(\sigma_+^{-\frac12}s^{-\frac12}z_N) & \mbox{ if } z_N > 0,\\
  f^\prime(\sigma_-^{-\frac12}s^{-\frac12}z_N) & \mbox{ if } z_N \le 0,
  \end{cases}
  \end{equation}
 where $f$ is given by \eqref{error function}. Once the limit function $S^*$ is uniquely determined, the compactness arguments again show that for every $n \in \mathbb N$:
 \begin{eqnarray}
  \{S^{\varepsilon}\}_{\varepsilon\le (n+2)^{-1}} \mbox{ converges to } S^* \mbox{ as } \varepsilon \to 0^+ \mbox{ uniformly on } \overline{Q_n}; &&\label{uniform convergence original sequence 2nd}\\
   \{\nabla S^{\varepsilon}\}_{\varepsilon\le (n+2)^{-1}} \mbox{ converges to } \nabla S^* \mbox{ as } \varepsilon\to 0^+ \mbox{ uniformly on } \overline{Q_{n\pm}}.&&\label{uniform convergence original sequence 2nd gradiet}
\end{eqnarray}
Namely, $S^*$ is regarded as the first-order approximation in $\varepsilon$  as $\varepsilon \to 0^+$ of $v^\varepsilon$.
In particular, by setting $s=1, \varepsilon =\sqrt{t}$ and $z=0$, we observe that 
\begin{equation}
\label{at q on the boundary}
S^\varepsilon(0,1) = t^{-\frac12}\!\!\left(u(0,t)-\frac{\sqrt{\sigma_+}}{\sqrt{\sigma_+}+\sqrt{\sigma_-}}\right) \mbox{ and }  S^*(0,1) =\frac{\sqrt{\sigma_+}\sqrt{\sigma_-}}{\sqrt{\pi}(\sqrt{\sigma_+}+\sqrt{\sigma_-})}(N\!-\!1)H(0),
\end{equation}
which gives  the conclusion \eqref{asymptotic formula in time} of Theorem \ref{th: asymptotic formula in time} from \eqref{uniform convergence original sequence 2nd}.

 \subsection{The uniform convergence on the interface}
 \label{uniform convergence}
 Let us prove  that the convergence  in \eqref{uniform convergence original sequence 2nd} at $(0,1)\in \mathbb R^{N}\times(0,+\infty)$ is uniform on $q \in \partial\Omega$, which suffices to prove the last conclusion of  Theorem \ref{th: asymptotic formula in time}. Suppose that this is not the case. Then there exist a number $\gamma > 0$, a sequence $\{\varepsilon_m\}$ with $\lim\limits_{m \to \infty}\varepsilon_m=0$, and a sequence $\{ q_m\} \subset \partial\Omega$ such that
 \begin{equation}
 \label{the contrary}
 |S^{\varepsilon_m}_{q_m}(0,1) - S^*_{q_m}(0,1)| \ge \gamma\ \mbox{ for every } m \in \mathbb N,
 \end{equation}
where $S^\varepsilon_{q_m}, S^*_{q_m}$ mean $S^{\varepsilon}, S^*$, whose $q \in \partial\Omega$ is replaced with $q_m\in\partial\Omega$, respectively. Recall that at their principal coordinate systems for $q_m\ (m \in \mathbb N)$, the origin $0$ corresponds to each $q_m$. Since $\partial\Omega$ is compact, it follows from the Bolzano-Weierstrass theorem that  $\{ q_m\}$ has a convergent subsequence. Thus, for brevity, we may assume that $\{ q_m\}$ itself converges to a point $q_*\in\partial\Omega$ as $m\to \infty$. Hence $H(q_m) \to H(q_*)$ as $m \to \infty$, since $\partial\Omega$ is of class $C^2$. By the second equality of \eqref{at q on the boundary}, we infer that
\begin{equation}
\label{one step to get a contradiction}
 S^*_{q_m}(0,1) \to  S^*_{q_*}(0,1)=\frac{\sqrt{\sigma_+}\sqrt{\sigma_-}}{\sqrt{\pi}(\sqrt{\sigma_+}+\sqrt{\sigma_-})}(N-1)H(q_*) \ \mbox{ as } m \to \infty,
\end{equation}
where $S^*_{q_*}$ means $S^*$ whose  $q \in \partial\Omega$ is replaced with $q_*\in\partial\Omega$.
Observe from \eqref{interior estimates for S varepsilon} that $\{ S^{\varepsilon_m}_{q_m}\}$ may also satisfy the interior estimates: for every $n, m \in \mathbb N$
\begin{equation*}
 \Vert S^{\varepsilon_m}_{q_m}\Vert_{C^{\alpha,\alpha/2}(\overline{Q_n})} + \Vert\nabla_{\hat{z}} S^{\varepsilon_m}_{q_m}\Vert_{C^{\alpha,\alpha/2}(\overline{Q_n})}+ \Vert \partial_{z_N}\! S^{\varepsilon_m}_{q_m}\Vert_{C^{\alpha,\alpha/2}(\overline{Q_{n+}})}+ \Vert \partial_{z_N}\! S^{\varepsilon_m}_{q_m}\Vert_{C^{\alpha,\alpha/2}(\overline{Q_{n-}})} \le C_{2,n}.\qquad
 \end{equation*}
Then it follows from the Arzel\`a-Ascoli theorem together with the Cantor diagonal process  that there exists a convergent subsequence of $\{ S^{\varepsilon_m}_{q_m}\}$. For brevity, let $\{ S^{\varepsilon_m}_{q_m}\}$ itself converge. In particular, since $q_m \to q_*$ as $m \to \infty$, as in the previous blow-up arguments, we may infer that as $m \to \infty$
 \begin{equation}
 \label{convergence to S_* at q_*}
 \{ S^{\varepsilon_m}_{q_m}\}\mbox{ converges to } S^*_{q_*}\ \mbox{ uniformly on every compact set in }\mathbb R^N\times(0,+\infty).
 \end{equation}
  Hence it follows from \eqref{convergence to S_* at q_*}  that
 \begin{equation}
 \label{the other step to get a contradiction}
 S^{\varepsilon_m}_{q_m}(0,1) \to S^*_{q_*}(0,1)=\frac{\sqrt{\sigma_+}\sqrt{\sigma_-}}{\sqrt{\pi}(\sqrt{\sigma_+}+\sqrt{\sigma_-})}(N-1)H(q_*) \ \mbox{ as } m \to \infty.
\end{equation}
This together with \eqref{one step to get a contradiction} yields a contradiction to \eqref{the contrary}.

\setcounter{equation}{0}
\setcounter{theorem}{0}

\section{Proofs of Theorems \ref{th:hyperplane} and \ref{th:sphere theorem}}
\label{section_Applications}

Let us prove Theorems  \ref{th:hyperplane} and \ref{th:sphere theorem}.

\noindent
{\it Proof of Theorem \ref{th:hyperplane}.}\ By Corollary \ref{constant mean curvature}, the mean curvature $H$ of $\partial\Omega$ must be constant. Since $\partial \Omega$ is a entire graph, by the divergence theorem $H$ must vanish.
Hence $\partial\Omega$ is an entire minimal graph over $\mathbb R^{N-1}$. If $N=2$, the mean curvature is simply a curvature, and hence $\partial\Omega$ must be a straight line. If \  $3\le N\le 8$, by the Bernstein theorem for the minimal surface equation (see \cite[Theorem 17.8, p.208]{G1984}), $\partial\Omega$ must be a hyperplane, and   if $\nabla \varphi$ is bounded with $ N \ge 3$, by Moser's theorem\cite[Corollary, p.591]{M1961} (see also \cite[Theorem 17.5, p.205]{G1984}), $\partial\Omega$ must be a hyperplane. \qed

\vskip 2ex

\noindent
{\it Proofs of Theorem \ref{th:sphere theorem}.}\ By Corollary \ref{constant mean curvature} together with Alexandrov's sphere theorem \cite[p.412]{Al1958}, $\Gamma$ must be a sphere. Since either the inside of $\Gamma$ or the outside of $\Gamma$ is included in one of the two sets, $\Omega$ and $\mathbb R^N\setminus\overline{\Omega}$, for instance let us consider the case where the inside of $\Gamma$ is included in $\mathbb R^N\setminus\overline{\Omega}$. Then, by taking into account  the initial data $\mathcal X_\Omega$ and the overdetermined condition \eqref{stationary isothermic surface with one C2 bounded part}, we infer from the uniqueness of the solution of the initial-boundary value problem for the heat equation $u_t=\sigma_-\Delta u$ that, for every $t >0$, $u$ is radially symmetric in $x$ with respect to the center of $\Gamma$ in the inside of $\Gamma$. Moreover,
by taking into account  the real analyticity in $x \in \Omega$  of the solution $u$ of the heat equation $u_t=\sigma_+\Delta u$ and  the transmission conditions \eqref{transmission conditions parabolic}, we see that, for every $t >0$, the radial symmetry  in $x$ of $u$ spreads into a connected component of $\Omega$ whose boundary includes $\Gamma$. 

Indeed, let $y\in \mathbb R^N$ be the center of $\Gamma$ and consider the partial differential operators $\mathcal R_{ij}$ given by
$$
\mathcal R_{ij}=(x_j-y_j)\partial_{x_i}-(x_i-y_i)\partial_{x_j}\ \mbox{ for } i\not=j.
$$
Then, since $\mathcal R_{ij}\Delta=\Delta \mathcal R_{ij}$, we see that each $\mathcal R_{ij}u$ also satisfies the same heat diffusion equation as $u$ does. Note that, for every $t >0$, $u$ is radially symmetric in $x$ with respect to $y$ in a domain $D$  if and only if $\mathcal R_{ij}u= 0$ in $D$ for every $i\not=j$. We observe that, for every $t >0$,  $\mathcal R_{ij}u= 0$ in the inside of the sphere $\Gamma$ for every $i\not=j$. By combining this fact with  the transmission conditions \eqref{transmission conditions parabolic} on $\Gamma  \times (0, +\infty)$, we see that the Cauchy data on $\Gamma \times (0,+\infty)$ for each  $\mathcal R_{ij}u$ must vanish. Hence it follows from Holmgren's uniqueness theorem that, for every $t >0$, each  $\mathcal R_{ij}u$ must vanish in a neighborhood of $\Gamma$ even in the outside of $\Gamma$. 
See for instance \cite[Theorem 4.2, p.250]{Mi1973} or \cite[Uniqueness Theorem, p.83]{J1982} for Holmgren's uniqueness theorem.
Furthermore,  the real analyticity in $x \in \Omega$  of the solution $u$ of the heat equation $u_t=\sigma_+\Delta u$ yields that, for every $t >0$, each $\mathcal R_{ij}u$ vanishes in a connected component of $\Omega$ whose boundary includes $\Gamma$. 

Suppose that  there is another component of $\partial\Omega$. Since $\partial\Omega$ is of class $C^0$, we may choose a component $\gamma$ which is the nearest to $\Gamma$. The spread radial symmetry in $x$ of $u$ and the overdetermined condition \eqref{stationary isothermic surface with one C2 bounded part} together with the real analyticity in $x \in \Omega$  of $u$ imply that $\gamma$ must be a sphere having the same center as that of $\Gamma$. Therefore $\gamma$ is of class $C^2$ and hence by  Theorem \ref{th: asymptotic formula in time} with the overdetermined condition \eqref{stationary isothermic surface with one C2 bounded part} the radius of $\gamma$ must be the same as that of $\Gamma$. This is a contradiction. The other cases can be dealt with similarly. \qed

\bigskip
\noindent{\Large\bf Acknowledgements.} The author would like to thank Professor Katsuyuki Ishii for the reference \cite{E1993} on the initial behavior of the solution $u$ of the Cauchy problem with initial data the characteristic function of a smooth open set for the heat equation $u_t=\Delta u$.


\end{document}